\documentclass{aimscleaned}

\usepackage{psfrag}
\usepackage{amsmath}
\usepackage{paralist}
\usepackage{graphics} 
\usepackage{epsfig} 
\usepackage{graphicx}  \usepackage{epstopdf}
\usepackage[colorlinks=true]{hyperref}
\hypersetup{urlcolor=blue, citecolor=red}

\def\currentvolume{X}

    \def\ppages{X--XX}
     \def\DOI{10.3934/xx.xx.xx.xx}

\theoremstyle{definition}

\newtheorem{remark}{Remark}

\usepackage{color}

\newcommand{\n}{\noindent}
\newcommand{\R}{\mathbb{R}}
\newcommand{\N}{\mathbb{N}}
\newcommand{\Lsp}{\mathbf{L}}
\newcommand{\PP}{\mathcal{P}}

\newcommand{\In}{\textsc{inc}}
\newcommand{\Out}{\textsc{out}}
\newcommand{\Start}{\textsc{start}}
\newcommand{\End}{\textsc{end}}
\newcommand{\Next}{\textsc{next}}
\newcommand{\rr}{\textsc{r}}
\newcommand{\jj}{\textsc{j}}
\newcommand{\dd}{\textsc{d}}

\newcommand{\Dx}{{\Delta x}}
\newcommand{\Dt}{{\Delta t}}
\newcommand{\rhomeno}{u_-}
\newcommand{\rhopiu}{u_+}

\title[A Destination-preserving Model for Wardrop Equilibria] 
      {A Destination-preserving Model for Simulating Wardrop Equilibria in Traffic Flow on Networks}

\author[E. Cristiani and F. S. Priuli]{}

\subjclass{Primary: 90B20; Secondary: 49N90.}
 \keywords{Traffic, networks, source-destination model, multi-path model, multi-population model, multi-commodity model, Wardrop equilibrium, Nash equilibrium.}

 \email{e.cristiani@iac.cnr.it}
 \email{f.priuli@iac.cnr.it}

\thanks{The research leading to these results has received funding from the European Union FP7 under grant No.\ 257462 HYCON2 Network of Excellence. The research was also supported by the collaboration with the private company ZEROPIU (Italy).}

\begin{document}
\maketitle

\centerline{\scshape Emiliano Cristiani and Fabio S. Priuli}
\medskip
{\footnotesize
 \centerline{Istituto per le Applicazioni del Calcolo ``M. Picone''}
   \centerline{Consiglio Nazionale delle Ricerche}
   \centerline{Via dei Taurini, 19}
   \centerline{I-00185 Rome, Italy}
} 

\bigskip

 \centerline{(Communicated by the associate editor name)}

\begin{abstract}
In this paper we propose a LWR-like model for traffic flow on networks which allows to track several groups of drivers, each of them being characterized only by their destination in the network. 
The path actually followed to reach the destination is not assigned \emph{a priori}, and can be chosen by the drivers during the journey, taking decisions at junctions.

The model is then used to describe three possible behaviors of drivers, associated to three different ways to solve the route choice problem: 1.\ Drivers ignore the presence of the other vehicles; 2.\ Drivers react to the current distribution of traffic, but they do not forecast what will happen at later times; 3.\ Drivers take into account the current and future distribution of vehicles. Notice that, in the latter case, we enter the field of differential games, and, if a solution exists, it likely represents a global equilibrium among drivers.

Numerical simulations highlight the differences between the three behaviors and offer insights into the existence of equilibria.
\end{abstract}

\section{Introduction}
In this paper we deal with a new LWR-like (i.e.\ macroscopic, differential, first-order) model for traffic flow on networks which allows one to track several groups of drivers, each of them being characterized only by a specific destination in the network.
The model allows one to simulate Wardrop and Nash equilibria in traffic flow.

\subsection{Background} In the past decades a great attention has been devoted to models which describe vehicular traffic flows via conservation laws. In such models, one studies the evolution of the density $\rho=\rho(t,x)$ of cars on a road, rather tracking each single car. The natural assumption that the total mass is conserved along the road leads to impose that $\rho$ obeys
\begin{equation}\label{eq:claw}
\partial_t \rho +\partial_x (\rho v(\rho)) = 0\,,
\qquad\qquad
\rho(0,x)=\rho_0(x)\,,
\end{equation}
for $(t,x)\in[0,+\infty[\times\R$ and for some initial distribution $\rho_0$. The function $v$ describes the dependence of the velocity of cars on their density, and it is typically chosen to be of the form $v(\rho)=V_*\left(1-{\rho\over \rho_*}\right)$ for some positive normalization constants $V_*, \rho_*$. This kind of first order models have been introduced by Lighthill and Whitham~\cite{LW} and by Richards~\cite{Richards}. 

However, in order to describe real situations where the cars move on a (typically very complex) network of roads, the simple model~\eqref{eq:claw} is not sufficient. This has motivated several authors to consider analogous equations on a \emph{network} $\mathcal{N}$, which is a directed graph whose nodes are called \emph{junctions} and whose arcs are called \emph{roads}.
The natural way to extend~\eqref{eq:claw} to a network $\mathcal{N}$ is to assume that the conservation law~\eqref{eq:claw} is separately satisfied on each road for all times $t>0$. Moreover, additional conditions have to be imposed at junctions, because in general the conservation of the mass alone is not sufficient to characterize a unique solution when two or more roads meet, even when the initial datum is piecewise constant. We refer the reader to the book by Garavello and Piccoli~\cite{piccolibook} for more details about the general ill-posedness of the problem at junctions. Multiple workarounds for such ill-posedness have been suggested in the literature: maximization of the fluxes across junctions (see~\cite{piccolibook}); buffer-like models where cars entering a junction with congested outgoing roads join a queue and exits in FIFO fashion (see~\cite{BressanKhai,BressanYu,GaravelloDCDS,GaravelloGoatin,herty2009NHM}); multi-path models which replace the junctions with suitable overlapping paths (see~\cite{bretti2014DCDS-S,briani2014sub}). In general, they all allow to determine a unique solution for the traffic evolution on the network, but the solution might be different.

Finally, traffic flow models on networks were extended to handle multiple groups of drivers, with different destinations (or objectives), sharing the network. The mathematical investigation on this subject has been very active in the past few years~\cite{BressanHan2,bretti2014DCDS-S,briani2014sub,Cong,garavello2005CMS}, both from the theoretical and the numerical point of view. The problem presents several difficulties and there is still no comprehensive theory capable to describe realistic networks.

\subsection{Goal}\label{sec:goal} 
The aim of this paper is twofold. On the one hand, we introduce a new differential model for multiple groups of drivers on a network. \textit{Each group is characterized only by a specific destination in the network}. The preferred path to be followed in order to reach the destination is not assigned \emph{a priori}, and can be modified during the journey, by making choices at junctions. Also, different drivers with the same destination (i.e.\ belonging to the same group) may find it convenient to use different paths at different times, because of, e.g.,  different traffic conditions. 

On the other hand, we exploit such a model to introduce different degrees of \emph{rationality} in drivers' choices, in the same spirit of~\cite{cristiani2014sub}. Namely, we couple the evolution equation with suitable control problems at junctions and we vary the amount of information that drivers can exploit in their decision procedure. In this way, we are able to describe both myopic contexts with low rationality involved, and high rationality contexts where the drivers are capable to find global equilibria/optima on the network. Numerical simulations, presented in section~\ref{sec:tests}, will show how the different degrees of rationality affect the final traffic flow. 

\medskip

In order to deal with multiple groups of drivers on $\mathcal{N}$ we proceed as follows. We assume that at each time and on each road,
the density of cars can be represented by
$$
\rho(t,x)=\sum_{d=1}^{N_\mathcal{D}} \rho_d(t,x)\,,
$$
where $N_\mathcal{D}$ is the number of destinations in $\mathcal{N}$, i.e.\ nodes of the network with no outgoing arcs, and $\{\rho_d\}_{d=1,\ldots,N_\mathcal{D}}$ are non-negative bounded functions representing the density of drivers at $(t,x)$ whose target is to eventually reach the $d$-th destination. Of course, assuming that the destination does not change during the travel and reasoning by linearity, it is natural to assume that each distribution satisfies
$$
\partial_t \rho_d +\partial_x (\rho_d\, v(\rho)) = 0\,,
\qquad\qquad
\rho_d(0,x)=\rho_{0,d}(x)
$$
outside junctions.

In order to allow drivers to change their path along the network during the journey, we assume to be given a family of functions $\Next_d(t,\textsc{j})$, for $d=1,\ldots,N_\mathcal{D}$, which describes, at each time $t$, what road to choose at junction $\textsc{j}$ to drive towards the $d$-th destination.
These functions can be seen as \textit{control parameters} to be either imposed \textit{a priori} by the network manager or to be chosen by the drivers as they travel through the network. They are only used at junctions because once a car has entered a certain road its path cannot be changed until next junction.

\medskip

The motivation to introduce a traffic flow model with the above characteristics comes from our interest in modeling different degrees of rationality in the drivers' behavior. Following the ideas presented in~\cite{cristiani2014sub}, in the context of pedestrian dynamics, we want to describe drivers which are capable to plan strategically the path toward their destination, making use of a certain amount of information. Offering different amounts of such information will lead to different strategies that can be interpreted as the result of different degrees of rationality in drivers' choices.

First of all, let us assume that each driver assigns a weight $w_\textsc{r}(t)$ to each road $\textsc{r}$ of the network, based on the information he/she has available. Then, at each time $t$ and each junction $\textsc{j}$, drivers aiming to the $d$-th destination node, denoted hereafter by $\dd_d$, choose $\Next_d(t,\textsc{j})$ as the road exiting from $\textsc{j}$ along the ``lighter'' path between $\textsc{j}$ and $\dd_d$. In other words, comparing the sum of weights of the various paths joining $\textsc{j}$ and $\dd_d$, $\Next_d(t,\textsc{j})$ will be the first road of the path realizing the minimum sum.

To fix the ideas, you can think that drivers want to get to their destination as soon as possible, and that the weight they assign to each road $\textsc{r}$ is simply the time necessary to drive through it, either ignoring or taking into account traffic conditions in $\textsc{r}$. In this case, at each time $t$ and each junction $\textsc{j}$, $\Next_d(t,\textsc{j})$ will be the outgoing road which ensures the shortest arrival time to the $d$-th destination.

We are now in the position to introduce more in detail the behaviors we can identify.

\begin{itemize}
\item {\bf basic behavior} (= no rationality): In this case, we assume that drivers have a very limited knowledge of the current status of the network and they can only choose their path on the basis of its geometry. Thus, each group will choose at start the functions $\Next_d$ which give the most convenient sequence of roads from their initial position to their destination \textit{as if no other cars are present}. In particular, the chosen path does never change during the evolution, because it is independent of the distribution of cars on $\mathcal{N}$.
\item {\bf rational behavior} (= mild rationality): In this case, we assume that drivers own some devices that reveals in real time the \textit{current} distribution of cars on each arc of the road, e.g.\ a smartphone, and that they can react to the modification of $\rho$ in the network by updating their controls $\Next_d$ accordingly. In this situation, at each time $\tau>0$, drivers can update the weights $w_\textsc{r}(\tau)$ of each road thanks to the knowledge of $\rho(\tau,\cdot)$ on $\mathcal{N}$ and select $\Next_d(\tau,\textsc{j})$ so as to optimize their path.
\item {\bf highly rational behavior} (= high rationality): In this case, drivers are not only informed of the \textit{current} distribution of cars on each arc of the network, but they can also \textit{forecast} accurately the evolution of such distribution and the long term effect of their own choices on the global evolution of traffic conditions. Here, the optimization problem depends on the whole distribution of $\rho$ in time and space, and thus it is fully coupled with the conservation law, like in the case of pedestrian flow models described via mean-field games~\cite{cristiani2014sub,dogbe2010MCM,LachaWolf,priuli2014sub}
\end{itemize}

Notice that, in the latter case, we enter the field of differential games, and, if a solution exists for the resulting coupled system, it likely represents a global equilibrium among drivers on the network. More precisely, if the weight assigned to each road simply corresponds to its travel time, we claim that the solution corresponding to this behavior satisfies the so-called Wardrop's first principle~\cite{Wardrop}:
\begin{quote}
The journey times on all the routes actually
used are equal, and less than those which
would be experienced by a single vehicle on
any unused route.
\end{quote} 
In particular, no individual driver can reduce his path cost by switching routes. 
Observe also that, since no unilateral change of strategy in some group of agents can lead to a reduction of the cost, Wardrop's equilibria can be seen as Nash equilibria among the drivers, see, e.g., \cite[Sect.\ 1]{carlier2012JMS}.

\medskip

\subsection{Comparison with the relevant literature} There are several novelties in our results compared to what is currently known in literature. First of all, the multi-population model that we present here is, to our knowledge, the first one allowing for time variable paths on the network. Indeed, early multi-population (or multi-class) models either do not consider networks at all, see, e.g.,~\cite{colombosbg} and references therein, or require populations to only use a fixed path to go from their source to their destination, see~\cite{garavello2005CMS, piccolibook} (cf.\ in particular~\cite[Def. 7.1.6, Thm. 7.2.1]{piccolibook}). The flexibility outlined above is not present in more recent models either: the multi-path model investigated in~\cite{bretti2014DCDS-S, briani2014sub} tracks different groups of drivers characterized by their path along the network, but paths cannot be modified at runtime; in~\cite{BressanHan3} each population can use multiple paths to get to their destination, but the choice has to be performed \emph{offline}, before the actual evolution starts, meaning that it cannot be modified afterwards.

A second innovative aspect concerns the type of optimization problems that we can study, thanks to the model we introduced. Several authors have presented traffic optimization for differential models on networks, but the point of view has been mostly the one of network managers. In other words, the overall flow on $\mathcal{N}$ is optimized w.r.t.\ some given criterion, but without accounting in any way the actual desires of drivers~\cite{cascone2007M3AS, cutolo2011, herty2006SIOPT, gugat2005JOTA, herty2003SISC}, which could be forced to take a much longer path before being able to reach their destinations. While this can be reasonable when dealing with very small networks, because in any case drivers can get back to their shortest path after exiting from the ``controlled'' part of the network, it seems unrealistic when dealing with large networks. On the other hand, some optimization results consider the desires of the drivers, and offer them the possibility to change path during the evolution, but only in the context of ``static'' models~\cite{Fisk,Nachtigall,OR1990}. Namely, on each road $\rr$ of the network, a cost $t\mapsto \mathcal{C}_\rr(t)$ is fixed and drivers can decide at each time which road to choose so to minimize the overall cost. However, in these models it is not clear how to pre-compute the functions $\mathcal{C}_\rr$, without considering the actual evolution of car densities on $\mathcal{N}$. For instance, it can be very hard to predict the appearance of congestions which might reduce the convenience of certain arcs. The approach we propose here, instead, allows to fill this gap, by coupling the optimization problem with the dynamics of the density.

A final aspect to consider is the connection between the solution to the coupled system for highly rational behaviors and Nash equilibria among drivers. The treatment of network flow as a differential game is not new. Some similar results are available in the literature, especially in the ``static'' models where a cost to pass through each road is given \textit{a priori} and no car dynamics is considered~\cite{Fisk,Hollander,Nachtigall,OR1990}. In our case the situation is much more complex because the costs depend on the distribution of cars. 
Results of different nature are contained in~\cite{carlier2008SIAM,carlier2012JMS}, where the distribution of cars on the network is seen as a transport plan between a mass of drivers concentrated at the origins and a mass of drivers concentrated at the destinations. In this case, equilibria on the network can be found as optimal transport plans, but the results are of stationary nature and it does not seem easy to include in the transport problem the presence of roads with different properties (capacity, maximal speed, etc.), which is instead almost straightforward in the models based on conservation laws.
 
More similar results are the ones presented in~\cite{BressanHan3,BressanKhai}. In the former paper, Nash equilibria among drivers are studied in detail from the theoretical point of view. Along the lines of the previous works~\cite{BressanHan1,BressanHan2}, the authors consider a model where players/drivers can choose their departure time (which acts as a control parameter) and the path they will follow on $\mathcal{N}$ so as to minimize a certain cost accounting for the arrival time and for the duration of the journey. It is proved the existence of a unique Nash equilibrium whenever the cost functional satisfies suitable regularity and monotonicity properties. The main difference between the problem studied here and the one in~\cite{BressanHan3} is again related to the possibility to change path during the journey: we allow drivers of each group to modify their initial choice at later time, and this is why our controls are defined only at the junctions, in the form of the functions $\Next_d$. As a drawback, we are not able to prove, at the moment, under which conditions the procedure we use to construct the Nash equilibrium does indeed converge to a solution of the problem. This is a very complex problem because it is quite common for differential games to have no Nash equilibria or, conversely, infinitely many equilibria (cf.~\cite{BrPr, CaCr, Pr}).

Finally, in~\cite{BressanKhai} a new multi-buffer model is introduced to handle the problems at junctions, producing a continuous (in $\Lsp^1$) semigroup of solutions for the traffic flow on the network. Such a model is very interesting and promising, but at the moment its numerical implementation is still under development, making difficult to use it in our context to compare which solutions are singled out when strategies with different rationality degrees are implemented.

\subsection{Paper organization} 
Section~\ref{sec:model} presents the destination-preserving model, while section~\ref{sec:rationality} explains how to introduce different degrees of rationality in the choices that drivers perform at junctions. Section~\ref{sec:coupling} focuses on how the two aspects can be coupled to describe rational traffic flows on a road network. In section~\ref{sec:numerics} we discretize the equations numerically, and in section~\ref{sec:tests} we present the result of some simulations which highlight the main differences between models with low rationality and models with a higher degree of rationality. Finally, section~\ref{sec:conclusions} presents conclusions and some open problems.

\section{The destination-preserving model}\label{sec:model}
We introduce here the model for the evolution of car densities on the network. 

\subsection{Preliminary notations and assumptions on the networks}\label{sec:notations}

In what follows, a network $\mathcal{N}$ will always be a directed graph consisting of a set of \emph{junctions} $\mathcal{J}$ (nodes) and a set of \emph{roads} $\mathcal{R}$ (arcs), i.e.,
$$
\mathcal{N} = \mathcal{J}\cup\mathcal{R}\,.
$$
We assume that for each junction $\textsc{j}\in\mathcal{J}$, there exist disjoint subsets
$$
\In(\textsc{j})\subset\mathcal{R}\,,\qquad\qquad
\Out(\textsc{j})\subset\mathcal{R}\,,
$$ 
representing, respectively, the incoming roads to $\textsc{j}$ and the outgoing roads from $\textsc{j}$. 
Among junctions, we distinguish two particular subsets consisting of \emph{origins} $\mathcal{O}$, which are the junctions $\textsc{j}$ such that $\In(\textsc{j})=\emptyset$, and \emph{destinations} $\mathcal{D}$, which are the junctions $\textsc{j}$ such that $\Out(\textsc{j})=\emptyset$. The junctions in $\mathcal{O}\cup\mathcal{D}$ can be considered as boundary points of $\mathcal{N}$.
We also denote by $N_\mathcal{R}$ the number of roads, by $N_\mathcal{J}$ the number of junctions, and by $N_\mathcal{D}$ the number of destinations.

We assume that each road $\rr\in\mathcal{R}$ can be seen as an interval $]a_\rr,b_\rr[\,\subset\R$. Given a road $\rr\in\mathcal{R}$, we will sometimes use $\Start(\rr)$ (resp. $\End(\rr)$) to indicate the junction corresponding to the infimum $a_\rr$ (resp. to the supremum $b_\rr$) of the interval $]a_\rr,b_\rr[$.

To avoid degeneracies, from now on we assume that in our network $\mathcal{N}$ the sets 
$
\mathcal{J}\,,
\mathcal{R}\,,
\mathcal{O}\,,
\mathcal{D}\,,
\mathcal{J}\setminus(\mathcal{O}\cup \mathcal{D})
$ 
are all non-empty, i.e., there are neither isolated nodes 
nor isolated roads. 
%
%
Moreover, it is convenient (but not strictly necessary) assuming that every destination is reachable from every origin.

\subsection{Basic ideas}\label{sec:basicideas}
First of all, let us assume to be given a family of functions $\Next_d\colon [0,+\infty[\times\big(\mathcal{J}\setminus\mathcal{D}\big)\to\mathcal{R}$, for $d=1,\ldots,N_\mathcal{D}$, such that, for every $t>0$,  
$$
\Next_d(t,\textsc{j})\in\Out(\textsc{j}).
$$ 
The role of $\Next_d(t,\textsc{j})$ is to prescribe, at each time $t$ and at each junction $\textsc{j}$, which road will be chosen next by drivers of the $d$-th group who are passing through the junction. Clearly, such functions can be seen as controls acting on the network and they can be either imposed by the traffic manager or chosen by the drivers themselves, depending on the situations that we want to model. In the former case, the functions can be arbitrary functions. 
In the latter case, the functions $\Next_d(t,\textsc{j})$ can be described in various ways, depending on the amount of information that we assume to be available to drivers. The actual details on how to define these functions in order to model drivers' preferences will be given in section~\ref{sec:rationality}.
It is useful to note here that functions $\Next_d$ act as bang-bang-like controls, meaning that, at any fixed time, they steer the groups of drivers to a \textit{single outgoing road}. As a consequence, the density $\rho_d$ of drivers belonging to the $d$-th group is not split in more than one outgoing road. This is not true instead for the total density $\rho$.


Furthermore, let us assume that  all functions $\{\Next_d\}_{d=1,\ldots,N_\mathcal{D}}$ are piecewise constant, i.e.\ that there exists a discretization of the interval $[0,+\infty[$, with step $\Delta\tau$ such that the functions $\Next_d$ are constant on each interval $[h\Delta\tau,(h+1)\Delta\tau[$, $h\in\N$. This will be needed later in the construction of the model.

\begin{figure}[th!]
\begin{center}
\begin{psfrags}
\psfrag{a1}{$a_1$} \psfrag{a2}{$a_2$}
\psfrag{b3}{$b_3$} \psfrag{b4}{$b_4$}
\psfrag{a1-d}{$b_1-\delta$} \psfrag{a2-d}{$b_2-\delta$}
\psfrag{a3+d}{$a_3+\delta$} \psfrag{a4+d}{$a_4+\delta$}
\psfrag{eq1}{Eq.\eqref{eq:clawd}}
\psfrag{eq2}{Eq.\eqref{eq:clawp}}
\includegraphics[width=0.8\textwidth]{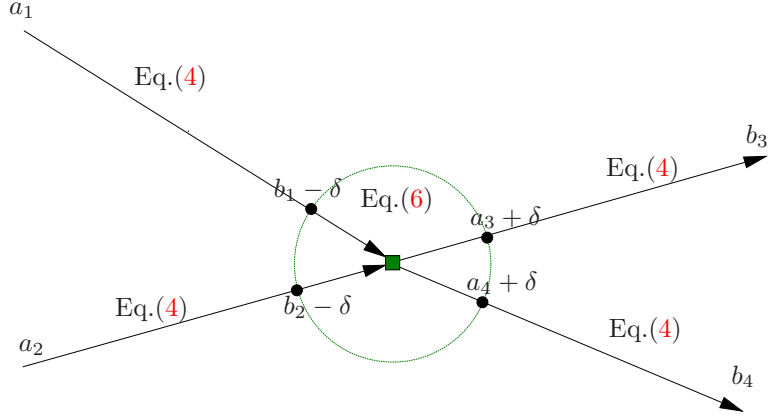}
\end{psfrags}
\end{center}
\caption{Separation between roads and junctions.}
\label{fig:interface}
\end{figure}

Finally, we have already mentioned that one of the most delicate parts of the model is the handling of the traffic flow across junctions. Here, we treat separately the flow away from the junctions and the flow in a neighborhood of each junction, see Fig.\ \ref{fig:interface}. 
In order to introduce rigorously such a separation, we fix a small parameter $\delta\in\left]0,\min_{\rr\in\mathcal{R}}{b_\rr-a_\rr\over 2}\right[$ and, for each road $\rr\in\mathcal{R}$, we define an open interval $I_\rr\subset\,]a_\rr,b_\rr[$ as follows
\begin{equation}\label{eq:inner_road}
I_\rr\doteq
\left\{
\begin{array}{ll}
]a_\rr,b_\rr-\delta[ &\mbox{ if }\Start(\rr)\in\mathcal{O}\,,\\&\\
]a_\rr+\delta,b_\rr[ &\mbox{ if }\End(\rr)\in\mathcal{D}\,,\\&\\
]a_\rr+\delta,b_\rr-\delta[ &\mbox{ otherwise}.
\end{array}
\right.
\end{equation}
In the following subsections, we introduce first the model in the \emph{interior of the roads} $\bigcup_{\rr\in\mathcal{R}} I_\rr$, and then the model in the \emph{proximity of each junction} $\textsc{j}\in\mathcal{J}\setminus(\mathcal{O}\cup \mathcal{D})$, i.e. in the sets 
\begin{equation}\label{eq:junct_neigh}
\mathcal{P}_{\textsc{j}}\doteq\bigcup_{i\in\In(\textsc{j})} [b_i-\delta,b_i[~~~~\cup \bigcup_{o\in\Out(\textsc{j})}\,]a_o,a_o+\delta[\,.
\end{equation}
To complete the definition of the model, suitable \emph{interface conditions} will be presented to ensure that the solutions obtained in the separate parts of the network do indeed give an admissible solution on the whole network $\mathcal{N}$.

Once the separate pieces are available, the problem on the whole network $\mathcal{N}$ is solved on the interval $[h\Delta\tau,(h+1)\Delta\tau]$, $h=0,1,\ldots$, and then the solution at $t=(h+1)\Delta\tau$ is used as initial datum for the problem on the interval $[(h+1)\Delta\tau,(h+2)\Delta\tau]$. 

\subsection{Interior of the roads} The modeling of traffic flow away from the junctions is done in the usual way. Each density distribution $\rho_d$, representing the drivers with destination $d$, evolves separately accordingly to the conservation law
\begin{equation}\label{eq:clawd}
\partial_t \rho_d +\partial_x (\rho_d\, v(\rho)) = 0\,,
\qquad
(t,x)\in \,]h\Delta\tau,(h+1)\Delta\tau[\,\times \bigcup_{\rr\in\mathcal{R}} I_\rr\,,
\end{equation}
where $\rho(t,x)\doteq \sum_{d=1}^{N_\mathcal{D}} \rho_d(t,x)$, complemented with suitable initial conditions $\rho_{d}(h\Delta\tau,\cdot)$ on the whole spatial domain and boundary conditions at $\mathcal O \cup \mathcal D$ at any time. 
We use here the classical Greenshield's model for the velocity, 
\begin{equation}\label{eq:fund_diag}
v(\rho)\doteq V_*\left(1-{\rho\over \rho_*}\right)\,,
\end{equation}
for fixed positive values $\rho_*,V_*$. The fact that the velocity depends on the whole distribution $\rho$ accounts for the fact that the velocity in the road is determined by the total amount of cars, independently of their destination. Notice that the model we are presenting here can be easily generalized to the case of different constants $\rho_{*,\rr},V_{*,\rr}$ for each road $\rr\in\mathcal{R}$, i.e.\ to the case of roads with different capacity and maximal speed. 


\begin{remark} We remark that, differently from other models presented in the literature, the density $\rho_d$ of drivers moving towards destination $d$ is defined in \emph{every road} of the network. It might well be that $\rho_d$ is identically zero in some roads. 
\end{remark}

\subsection{Junctions} \label{sec:junctions}
Let us consider a generic junction $\textsc{j}\in\mathcal{J}\setminus(\mathcal{O}\cup \mathcal{D})$. We observe that $\textsc{j}$, together with the corresponding neighborhood $\mathcal{P}_{\textsc{j}}$ defined in~\eqref{eq:junct_neigh}, can be seen as a simplified network, denoted hereafter by $\mathcal{N}_{\textsc{j}}$, consisting of a single junction $\textsc{j}$ and $|\In(\textsc{j})| + |\Out(\textsc{j})|$ roads. 

In our framework, the discretized functions $\{\Next_d\}_d$ determine a unique path on the network $\mathcal{N}_{\textsc{j}}$ for each group of drivers on each incoming road, consisting of the incoming road itself and of the outgoing road given by $\Next_d$.  Moreover, such a path remains fixed on the whole time interval $[h\Delta\tau,(h+1)\Delta\tau[$, $h\in\N$, because $\{\Next_d\}_d$ are constant functions in such interval.
Therefore, in order to model the evolution of the densities on $\mathcal{N}_{\textsc{j}}$, we can apply one of the several existing approaches like, e.g., the source-destination model~\cite{garavello2005CMS}, the model proposed in~\cite{BressanKhai}, or the multi-path model~\cite{bretti2014DCDS-S, briani2014sub}. 

In the following we shall employ the multi-path model since it is by far the simplest one among the cited ones, not requiring any separate procedure to compute the flux through the junction (e.g., the maximization of the flux). However, in order to apply such a model we need to define suitably what ``populations'' (in the sense of \cite{bretti2014DCDS-S, briani2014sub}) are in our context: indeed, the multi-path model tracks the evolution of different ``populations'' of drivers, \textit{characterized by their path} along the network. 

In our framework, we can identify up to $|\In(\textsc{j})|\times N_\mathcal{D}$ paths (=populations) on the subnetwork $\mathcal{N}_{\textsc{j}}$, corresponding to all admissible combinations incoming-outgoing roads at $\jj$ (here, ``admissible'' means that outgoing roads for drivers with destination $\dd_d$ are given by $\Next_d$).

We are finally in the position to formulate our model on $\mathcal{N}_{\textsc{j}}$ on each interval of times $[h\Delta\tau,(h+1)\Delta\tau[$. For any pair $(i,d)\in \In(\textsc{j})\times \{1,\ldots,N_\mathcal D\}$, we consider a new unknown $\mu^{\textsc{j}}_{(i,d)}$ which represents the density of a single population of drivers. The new unknown is defined on the path $\PP^{\jj}_{(i,d)}$ on $\mathcal{N}_{\textsc{j}}$, which consists of the concatenation of the last part of the incoming road $i$ and the first part of the outgoing road $\Next_d(t,\textsc{j})$, i.e., 
$$
\PP^\jj_{(i,d)}\doteq [b_i-\delta,b_i]\cup[a_o,a_o+\delta], \qquad o=\Next_d(t,\textsc{j}).
$$
Note that two or more paths can share some portion of roads. 
Following \cite{bretti2014DCDS-S, briani2014sub}, on the subnetwork $\mathcal{N}_{\textsc{j}}$ we have to solve the following system of $|\In(\textsc{j})|\times N_\mathcal{D}$ conservation laws with space-dependent and discontinuous flux 
\begin{equation}\label{eq:clawp}
\partial_t \mu^{\textsc{j}}_{(i,d)} +\partial_x \big(\mu^{\textsc{j}}_{(i,d)}\, v(\mu^\jj)\big) = 0\,,
\quad
(t,x)\in \,]h\Delta\tau,(h+1)\Delta\tau[\,\times \PP^\jj_{(i,d)}\,,
\quad \forall (i,d)
\end{equation}
where $\mu^\jj$ is the sum of all the densities defined on $\PP^\jj_{(i,d)}$ living at $x$ at the same time, i.e.
$$
\mu^\jj(t,x)\doteq 
\sum_{i\in\In(\textsc{j})}
\sum_{d=1}^{N_\mathcal{D}} \mu^{\textsc{j}}_{(i,d)}(t,x)\,, \qquad
(t,x)\in \,]h\Delta\tau,(h+1)\Delta\tau[\,\times \PP^\jj_{(i,d)}\,,
$$ 
and the velocity is still given by \eqref{eq:fund_diag}. Note that in this formulation the junction (apparently) disappears since each path is considered as a \textit{single uninterrupted one-dimensional domain}. Actually, the junction is hidden in the function $\mu^\jj$, which couples the equations and has a discontinuity at the junction, see \cite{bretti2014DCDS-S, briani2014sub} for details.

The system \eqref{eq:clawp} is complemented by initial conditions $\mu^{\textsc{j}}_{(i,d)}$ at time $t=h\Delta\tau$, on each path $\PP^\jj_{(i,d)}$. If we are given only the initial conditions for the densities $\rho_d$'s, we need to distribute such a value in the outgoing roads among the populations $(i,d)$ with $i\in\In(\textsc{j})$, i.e.\ we need to choose coefficients 
\begin{equation}\label{eq:pop_percent}
\Lambda^{\textsc{j}}_{(i,d)}\in [0,1]
\qquad
\mbox{ such that }
\sum_{i\in\In(\textsc{j})}\Lambda^{\textsc{j}}_{(i,d)}=1
\end{equation}
and define $\mu^{\textsc{j}}_{(i,d)}=\Lambda^{\textsc{j}}_{(i,d)}\rho_d$ in the outgoing part of the path $\PP^\jj_{(i,d)}$. A possible choice for the $\mu^{\textsc{j}}_{(i,d)}(h\Delta\tau,\cdot)$'s is the following: for each pair $(i,d)$, we define
$$
\Lambda^{\textsc{j}}_{(i,d)}\doteq 
{\displaystyle\int_{b_i-\delta}^{b_i}\mu^{\textsc{j}}_{(i,d)}(h\Delta\tau,x)\,dx
\over 
\displaystyle \sum_{i\in\In(\textsc{j})} \int_{b_i-\delta}^{b_i}\mu^{\textsc{j}}_{(i,d)}(h\Delta\tau,x)\,dx}
$$ 
and we set
$$
\mu^{\textsc{j}}_{(i,d)}(h\Delta\tau,x)\doteq \left\{
\begin{array}{ll}
\rho_{d}(h\Delta\tau,x) & \mbox{ if } x\in [b_i-\delta,b_i]\,,\\
\Lambda^{\textsc{j}}_{(i,d)}\,\rho_{d}(h\Delta\tau,x) & \mbox{ otherwise.}
\end{array}
\right.
$$

\subsection{Mass conservation at interfaces} It remains to specify how the two models are matched at the boundaries of the subnetwork $\mathcal{N}_\textsc{j}$. Luckily, the conservation of mass across the interfaces $x=b_i-\delta$ and $x=a_o+\delta$
is sufficient to determine uniquely the solution at the interfaces. 
Indeed, if we are given the solutions $\{\mu^{\textsc{j}}_{(i,d)}\}_{(i,d)}$ at every time $t>0$ in the subnetwork  $\mathcal{N}_\textsc{j}$, we determine the boundary conditions for the $\{\rho_d\}_d$ in the interior of the roads as follows:
\begin{itemize} 
\item for all $i\in\In(\textsc{j})$ and all $t>0$, $\rho_d(t,b_i-\delta) = \mu^{\textsc{j}}_{(i,d)}(t,b_i-\delta)$;
\item for all $o\in\Out(\textsc{j})$ and all $t>0$, 
$$
\rho_d(t,a_o+\delta) =\left\{
\begin{array}{ll}
\displaystyle\sum_{i\in\In(\textsc{j})} \mu^{\textsc{j}}_{(i,d)}(t,a_o+\delta) & \mbox{ if } o=\Next_d(t,\textsc{j})\,,\\
0 & \mbox{ otherwise.}
\end{array}
\right.
$$
\end{itemize}

\n Similarly, if we are given the solutions $\{\rho_d\}_d$ at every time $t>0$ in the domain $\bigcup_{\rr\in\mathcal{R}} I_\rr$, we determine the boundary conditions for the $\{\mu^{\textsc{j}}_{(i,d)}\}_{(i,d)}$ in each subnetwork $\mathcal{N}_\textsc{j}$ as follows:
\begin{itemize} 
\item for all $i\in\In(\textsc{j})$ and all $t>0$, $\mu^{\textsc{j}}_{(i,d)}(t,b_i-\delta)=\rho_d(t,b_i-\delta)$;
\item for all $o\in\Out(\textsc{j})$ and all $t>0$ choose again $\Lambda^{\textsc{j}}_{(i,d)}$ with the properties in~\eqref{eq:pop_percent} and define
$\mu^{\textsc{j}}_{(i,d)}(t,a_o+\delta) = \Lambda^{\textsc{j}}_{(i,d)}\, \rho_d(t,a_o+\delta)$.
\end{itemize}

\section{Modeling drivers' choice processes via path optimization} \label{sec:rationality}
In this section we describe the procedure for determining the control functions $\Next_d$. As mentioned in the introduction, our goal is to describe different decision processes to model different degrees of drivers' rationality. The common framework for all levels of rationality is that decisions are taken in order to optimize some performance criterion that drivers have. \emph{Consistently with the fact that drivers' choices are made only at junctions, the optimization on $\mathcal{N}$ will be of discrete type}.

Let us assume that it is given any total distribution of cars along the roads $\mathcal{R}$ for all times $[t_0,+\infty[\,$, i.e. any function 
$$
\omega\colon[t_0,+\infty[\,\times\bigcup_{\rr\in\mathcal{R}} ]a_\rr,b_\rr[\to [0,\rho_*],
$$
and that drivers with any destination assign some common weights $w_\rr[\omega]$ to each road $\rr\in\mathcal{R}$, which in general will depend on the values of $\omega$ in the whole road $\rr$ for all times. A natural choice to define such weights is to assume that they are computed as follows:\\
{\bf a.} since the drivers who start at time $t_0$ from the beginning of the road $\Start(\rr)$ expect to travel along the road $\rr$ according to the microscopic dynamics
\begin{equation}\label{eq:ode}
\dot x(s) = v(\omega(s,x(s))),\qquad\qquad x(t_0)=a_\rr,
\end{equation}
$v$ being the velocity field in~\eqref{eq:fund_diag}, we define the \emph{arrival time} at $\End(\rr)$ implicitly as the value $t_f=t_f[\omega](t_0)\in[t_0,+\infty]$ such that the solution of the previous ODE satisfies $x(t_f)=b_\rr$;\\
{\bf b.} we now define the weight
\begin{equation}\label{eq:weight}
w_\rr[\omega](t_0):=\int_{t_0}^{t_f[\omega](t_0)} \ell_\rr(s\,;\omega(s,x(s)))\,ds,
\end{equation}
for some suitable nonnegative running costs $\ell_\rr$. 
To fix the ideas, you can think to the special case of weights which correspond to the travel times through each road, namely $\ell_\rr\equiv 1$.

Notice that the generalization to the case of weights that also depend on the destination $\dd_d$ of each group of drivers is straightforward and left to the interested reader: it is a matter of replacing the $w_\rr$ with suitable new weights $w_{d,\rr}$, and to adapt accordingly the definition of the functions $V_d$ below.

Given a junction $\textsc{j}\in\mathcal{J}$ and a destination $\dd\in\mathcal{D}$, we call \emph{path joining $\textsc{j}$ with $\dd$} any sequence of arcs $\rr_1,\ldots,\rr_N \in \mathcal{R}$ such that $\Start(\rr_1)=\textsc{j}$, $\End(\rr_N)=\dd$ and $\End(\rr_k)=\Start(\rr_{k+1})$ for $k=1,\ldots,N-1$. Moreover, in such a case we call \emph{weight of the path} the sum $\sum_{k=1}^N w_{\rr_k}$. Of course, whenever the weights $w_\rr$ vary in time, due to a distribution $\omega$ that varies in time, then also the weight of the path will change in time.

By using the weights of paths, it is natural to define a \emph{value function} for the group of drivers aiming at the destination $\dd_d$ as the function $V_d\colon [0,+\infty[\times\mathcal{J}\to\R\cup\{+\infty\}$ such that $V_d(t,\jj)$ gives the minimum weight of the path among all paths joining $\textsc{j}$ with $\dd_d$ at time $t$. Functions $\Next_d(t,\textsc{j})$ can be then defined by choosing the outgoing road that belongs to the path realizing $V_d(t,\textsc{j})$.

To make things more precise, we distinguish three different processes that drivers can use to define their weights $w_\rr$ during the evolution, based on the different amounts of information that they can exploit. 

\subsection{Basic behavior} \label{sec:base}
In case of \emph{basic behavior} we simply assume that each driver chooses his/her optimal next road at junctions ignoring the presence of other drivers in the network, i.e.\ without considering the effect of the distribution $\rho$ given on $\mathcal{N}$. This corresponds to compute the weights as $w_\rr[\rho\equiv 0]$. 
In this case, the weights are the same for all times, and therefore the optimal path is time-independent as well. 
By means of the well known Dynamic Programming Principle, one can characterize the value function $V_d$ as the solution to the system
\begin{equation}\label{eq:HJBbasic}
V_d(\textsc{j})=\min_{\rr\in \Out(\textsc{j})} \big\{V_d(\End(\rr))+ w_\rr[0]\big\}\,,
\qquad \forall \jj\,,
\end{equation}
with boundary conditions on $\mathcal{D}$ given by $V_d(\dd_d)=0$ and $V_d(\dd_e)=+\infty$ for $e\in\{1,\ldots,N_\mathcal{D}\}\setminus\{d\}$. Here, $w_\rr[0]$ is as in~\eqref{eq:weight} with $t_0=0$, $\omega\equiv 0$, and thus $t_f=\,{b_\rr-a_\rr\over V_*}$.
Once $V_d$ has been computed, the optimal control is easily found as 
$$
\Next_d(t,\textsc{j})\doteq \arg\!\!\!\min_{\rr\in \Out(\textsc{j})} \big\{V_d(\End(\rr))+ w_\rr[0]\big\}\,,
$$
which is thus constant in time.
%

\subsection{Rational behavior}\label{sec:rat}
In case of \emph{rational behavior} we assume that at any fixed time $\tau\geq 0$ each driver is aware of the distribution $\rho(\tau,\cdot)$ in the whole network $\mathcal{N}$ and uses this information to select his/her own optimal path towards the destination. Then, treating $\tau$ as a fixed parameter, we define the value function $V_{d,\tau}$ as the solution to the system
\begin{equation}\label{eq:HJBrational}
V_{d,\tau}(\textsc{j})=\min_{\rr\in \Out(\textsc{j})} \big\{V_{d,\tau}(\End(\rr))+ w_\rr[\rho(\tau,\cdot)](\tau)\big\}\,,
\qquad\forall \jj\,,
\end{equation}
with boundary conditions on $\mathcal{D}$ given by $V_{d,\tau}(\dd_d)=0$ and $V_{d,\tau}(\dd_e)=+\infty$ for $e\in\{1,\ldots,N_\mathcal{D}\}\setminus\{d\}$. Here, $w_\rr[\rho(\tau,\cdot)](\tau)$ is as in~\eqref{eq:weight} with $t_0=\tau$ and $\omega=\rho(\tau,\cdot)$, so that $w_\rr[\rho(\tau,\cdot)](\tau)=\int_\tau^{t_f}\ell_\rr(s\,;\rho(\tau,x(s)))\,ds$ and $x(\cdot)$ denotes the (Carath\'eodory) solution to~\eqref{eq:ode}\footnote{Observe that if any portion of the road is fully congested ($\rho=1$), then $t_f=+\infty$ because cars entering at $a_\rr$ have no way to reach $b_\rr$. On the other hand, whenever $\rho<1$ in the whole road, $t_f<+\infty$ and $v>0$, and this implies that a Carath\'eodory solution to~\eqref{eq:ode} exists for positive times because no trajectory can remain trapped in a switching point of the vector field $v$.}.
Once $V_{d,\tau}$ has been computed, the optimal control is easily found as 
$$
\Next_d(\tau,\textsc{j})\doteq \arg\!\!\!\min_{\rr\in \Out(\textsc{j})} \big\{V_{d,\tau}(\End(\rr))+ w_\rr[\rho(\tau,\cdot)](\tau)\big\}\,.
$$
Repeating the construction for every $\tau>0$, we obtain the desired function $\Next_d$ at any time $t=\tau$.

\subsection{Highly rational behavior}\label{sec:high}
In case of \emph{highly rational behavior} we assume that drivers can exploit the complete knowledge of the distribution of cars $\rho$ in the whole network $\mathcal{N}$ at any time. Then, we define the value function $V_d$ as the solution to the system
\begin{equation}\label{eq:HJBhrational}
V_d(t,\textsc{j})=\min_{\rr\in \Out(\textsc{j})} \big\{V_d(t+\Upsilon_\rr[\rho(\cdot,\cdot)](t),\End(\rr))+ w_\rr[\rho(\cdot,\cdot)](t)\big\}\,,
\quad
\forall t\geq 0\,,~~ \forall \jj\,,
\end{equation}
with boundary conditions on $\mathcal{D}$, for all $t\geq 0$, given by $V_d(t,\dd_d)=0$ and $V_d(t,\dd_e)=+\infty$ for $e\in\{1,\ldots,N_\mathcal{D}\}\setminus\{d\}$. In~\eqref{eq:HJBhrational} the notation $w_\rr[\rho(\cdot,\cdot)](t)$ is used exactly like in~\eqref{eq:weight} with $t_0=t$, and $\Upsilon_\rr[\rho(\cdot,\cdot)](t)
$ denotes the time needed to drive through road $\rr$, when the total density $\rho$ is accounted along the arc (equivalently $\Upsilon_\rr[\rho(\cdot,\cdot)](t)$ can be thought as given by~\eqref{eq:weight} with $t_0=t$ and $\ell_\rr\equiv 1$). Once $V_d$ has been computed, the optimal control is easily found as 
\begin{equation}\label{eq:HJBhrational_next}
\Next_d(t,\textsc{j})\doteq \arg\!\!\!\min_{\rr\in \Out(\textsc{j})} \big\{V_d(t+\Upsilon_\rr[\rho(\cdot,\cdot)](t),\End(\rr))+ w_\rr[\rho(\cdot,\cdot)](t)\big\}\,.
\end{equation}

The difficulty in this case lies on the fact that the equation used to define the value function $V_d(t,\textsc{j})$ is \emph{fully coupled} with the model for the evolution of the density, since the functions $V_d$ affect the distribution $\rho$ via $\Next_d$, and, in turn, the distribution $\rho$ affects $V_d$. In particular, availability of information on the distribution $\rho$ in the whole time-space domain means that drivers do forecast the other drivers' decisions as well as the effect of their own choices on the drivers with different destinations. This will lead us naturally to Nash equilibria among the drivers as if they were players of a differential game.

%
%
%
%
%
\section{Coupling the evolution model with the path optimization}\label{sec:coupling}
In section~\ref{sec:model}, we have presented how to construct a time-discrete solution of~\eqref{eq:claw} on the network $\mathcal{N}$, whenever piecewise constant functions $\Next_d$ are given. On the other hand, in section~\ref{sec:rationality}, we have presented several ways to define the functions $\Next_d$, whenever a distribution of cars $\rho$ is given in $[0,+\infty[\,\times\mathcal{N}$. Here, we want to show how to combine the two constructions so as to describe different behaviors of the drivers in the network. In particular, we focus our attention on the coupling between the conservation laws~\eqref{eq:clawd},\eqref{eq:clawp} and the Hamilton-Jacobi equation \eqref{eq:HJBbasic} or \eqref{eq:HJBrational} or \eqref{eq:HJBhrational}, and on the discretization-in-time procedure of the functions $\Next_d$, which is needed to apply the model described in section~\ref{sec:model}.

We face different situations depending on the behavior we are trying to model.
\begin{itemize}
\item \emph{Basic behavior}: In this case, the functions $\Next_d$ are constant in time and their construction only depends on the ``geometry'' of the network; as such, there is no problem in coupling the two construction, since the discretization of each function $\Next_d$ is the function itself.
\item \emph{Rational behavior}: In this case, each function $\Next_d$ at time $t$ only depends on the distribution of cars $\rho(t,\cdot)$ at that same time. Therefore, we can procede iteratively: if we are given $\rho$ and $\Next_d$ on the interval $[0,h\Delta\tau]$, we can use $\rho(h\Delta\tau,\cdot)$ to construct $\Next_d(h\Delta\tau,\cdot)$, for all $d$; and then keep this function constant in $[h\Delta\tau,(h+1)\Delta\tau]$, so to be able to solve~\eqref{eq:clawd},\eqref{eq:clawp} up to time $(h+1)\Delta\tau$.
\item \emph{Highly rational behavior}: This is the most delicate case, because the functions  $\Next_d$ depend on the whole distribution $\rho(\cdot,\cdot)$ on $]0,+\infty[\,\times\mathcal{N}$, and, in turn, $\rho$ is affected by any change of the functions $\Next_d$. In order to construct a solution to the coupled system~\eqref{eq:clawd},\eqref{eq:clawp},\eqref{eq:HJBhrational},\eqref{eq:HJBhrational_next}, we shall look for fixed points of the following operator $\Xi$: given a distribution $\rho$ on $]0,+\infty[\,\times\mathcal{N}$, we use such $\rho$ to construct the functions $\{\Next_d\}_d$ as described in section~\ref{sec:high}. Then, we discretize the resulting functions by setting
$$
\Next^{h}_d(t,\textsc{j})\equiv \Next_d(h\Delta\tau,\textsc{j})
\qquad
\forall\,t\in[h\Delta\tau,(h+1)\Delta\tau]
$$
and we use such piecewise constant functions to construct a new distribution on $]0,+\infty[\,\times\mathcal{N}$ that we denote with $\Xi(\rho)$. At the moment, we are not able to say under which conditions the operator $\rho\mapsto\Xi(\rho)$ has fixed points. In next sections we start investigating numerically the behavior of the operator $\Xi$, by studying the behavior of iterated applications $\Xi(\rho)$, $\Xi^2(\rho)$, and so on.
\end{itemize}

%
%
%
%
%
\section{Numerical approximation}\label{sec:numerics}
Equations \eqref{eq:clawd} and \eqref{eq:clawp} have the same structure and can be approximated by the same numerical scheme. Let us unify the numerical handling of the two equations introducing a generic system of the form
\begin{equation}\label{eq:clawu}
\partial_t u_\alpha + \partial_x \big( u_\alpha v(u)\big),\qquad (t,x)\in]0,+\infty[\times Q_\alpha, \qquad \alpha=1,\ldots,N_\alpha\,,
\end{equation}
where:
\begin{itemize}
\item 
$\{Q_\alpha\}_{\alpha=1,\ldots N_\alpha}$ are $N_\alpha>0$ one-dimensional domains, possibly coinciding (as in \eqref{eq:clawd}) or having some parts in common (as in \eqref{eq:clawp});
\item 
$u_\alpha$ is the density of cars moving along $Q_\alpha$;
\item
$u$ is the sum of all the densities living at some point of $Q_\alpha$ at the same time, i.e. 
$$
u(t,x):=\sum_\alpha u_\alpha(t,x).
$$
\end{itemize}
Equation \eqref{eq:clawu} is complemented with suitable initial and boundary conditions.

In order to employ the Godunov-based discretization proposed in \cite{bretti2014DCDS-S,briani2014sub}, it is convenient to introduce the flux function 
$$
f(u):=uv(u),
$$
with $v$ as in~\eqref{eq:fund_diag}, and rewrite \eqref{eq:clawu} as
\begin{equation}\label{eq:clawuf}
\partial_t u_\alpha + \partial_x \left( \frac{u_\alpha}{u} f(u)\right),\qquad (t,x)\in]0,+\infty[\times Q_\alpha, \quad \alpha=1,\ldots,N_\alpha,
\end{equation}
setting $\frac{u_\alpha}{u}=0$ if $u=0$ (and then $u_\alpha=0$ $\forall \alpha$) to avoid singularities.

\medskip

We define a numerical grid in $[0,+\infty[\times Q_\alpha$ with space step $\Dx$ and time step $\Dt$.
We denote by $x_k:=k\Delta x$, $k\in\mathbb Z$, the center of the $k$-th space cell along $Q_\alpha$, and by $t^n:=n\Delta t$, $n\in\mathbb N$, the center of the $n$-th time cell. 
We also denote by $u^{n,k}_\alpha$ the approximate density $u_\alpha(x_k,t^n)$ and we naturally define
\begin{equation}\label{def:omega_approx}
u^{n,k}:=\sum_{\alpha}u^{n,k}_\alpha.
\end{equation}
Equation (\ref{eq:clawuf}) is discretized by means of the following Godunov-type scheme \cite{bretti2014DCDS-S,briani2014sub}, which reads, at any internal cell $k$, as
\begin{equation}\label{schema}
u^{k,n+1}_\alpha = u^{k,n}_\alpha-
\frac{\Dt}{\Dx}\left(\frac{u^{k,n}_\alpha}{u^{k,n}} \ G(u^{k,n},u^{k+1,n})- 
\frac{u^{k-1,n}_\alpha}{u^{k-1,n}}\ G(u^{k-1,n},u^{k,n})
\right)
\end{equation}
for $n\geq 0$ and $\alpha=1,\ldots,N_\alpha$, where $G$ is the classical Godunov numerical flux defined, as usual, as
\begin{equation}\label{GodunovFlux}
G(\rhomeno,\rhopiu):=
\left\{
\begin{array}{ll}
\min_{z\in [\rhomeno,\rhopiu]} f(z) & \textrm{if } \rhomeno\leq \rhopiu, \\
\max_{z\in [\rhopiu,\rhomeno]} f(z) & \textrm{if } \rhomeno\geq \rhopiu. \\
\end{array}
\right.
\end{equation}

The scheme \eqref{schema} has been proven to hold some nice properties. In particular, no special management of the junctions (i.e., the points where two or more $Q_\alpha$'s meet) is needed since the scheme selects \textit{automatically} a solution at junctions that maximizes the flow  along each path $Q_\alpha$ (user optimum). The scheme does not compute in general the maximal flow that could possibly be transferred over the node (global optimum), as it happens in more standard approaches \cite{piccolibook}. Moreover, when the demand of the incoming roads is larger than the supply of the outgoing roads (i.e. queues are formed behind the junction), the scheme equidistributes the incoming flux among the incoming roads, giving to the incoming roads the same priority. See \cite{bretti2014DCDS-S,briani2014sub} for more details.

\medskip

Regarding the optimization problem, we simply apply a fixed-point algorithm starting from the following initial guess:
$$
\left\{
\begin{array}{ll}
V^{\textup{guess}}_d(\jj)=0, & \jj= \dd_d \\
V^{\textup{guess}}_d(\jj)=+\infty, & \jj\neq \dd_d \\
\end{array}
\right.
$$
and then iterating the computation \eqref{eq:HJBbasic} or \eqref{eq:HJBrational} or \eqref{eq:HJBhrational} for all $\jj\in\mathcal J\backslash\mathcal D$ until convergence is reached.

The case of the highly rational behavior clearly has some additional complications, as already explained in section~\ref{sec:coupling}. First of all, since we are forced to set a final time $T$ for the simulation, the value $t+\Upsilon_\rr(t)$ appearing in \eqref{eq:HJBhrational} can be larger than $T$ for some $t$. If this is the case, the value function is no longer defined, therefore neither is the function $\Next_d$. We overcome this problem stopping the inflow of new cars in such a way that the network empties before a certain time, thus making not influential the fact that the value function is undefined at later times.

Second, it is not guaranteed that the iterations between the forward-in-time equations \eqref{eq:clawd},\eqref{eq:clawp} and the backward-in-time equation \eqref{eq:HJBhrational} converge to a solution. Numerical evidence shows that this is not always true: in some cases the algorithm oscillates between two solutions, possibly two Wardrop equilibria for the system. Cf.\ on this point the results in \cite{CaCr} for the numerical approximation of the differential games studied in~\cite{BrPr, Pr}.

We also remark that the initial guess used to trigger the fixed-point iterations for the operator $\Xi$ (see section \ref{sec:coupling}) affects the final results. In this paper we have chosen the density corresponding to the basic behavior as first guess. The density corresponding to the rational behavior can be also used instead.
%
%
%
%
%
\section{Numerical tests}\label{sec:tests}
\begin{figure}[h!]
\begin{center}
\begin{psfrags}
\psfrag{r1}{$\!\!\!\rr_1$}\psfrag{r2}{$\rr_2$}\psfrag{r3}{$\rr_3$}\psfrag{r4}{$\rr_4$}
\psfrag{r5}{$\rr_5$}\psfrag{r6}{$\rr_6$}\psfrag{r7}{$\rr_7$}\psfrag{r8}{$\rr_8$}
\psfrag{j1}{$\jj_1$}\psfrag{j2}{$\jj_2$}\psfrag{j3}{$\jj_3$}\psfrag{j4}{$\jj_4$}
\psfrag{j5}{$\jj_5$}\psfrag{j6}{$\!\!\!\jj_6$}\psfrag{j7}{$\jj_7$=$\dd_1$}\psfrag{j8}{$\jj_8$=$\dd_2$}
\includegraphics[width=0.8\textwidth]{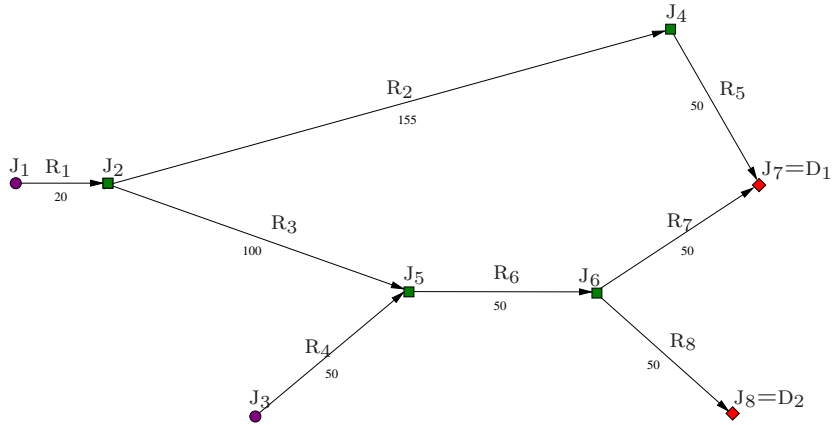}
\end{psfrags}
\end{center}
\caption{The network considered for the numerical tests (2 origins, 2 destinations, 8 roads, 8 junctions).}
\label{fig:network}
\end{figure}

In this section we present four numerical tests to prove the feasibility of our model and to show the differences among the three behaviors described in section \ref{sec:goal}.
We choose $V_*=\rho_*=1$ in \eqref{eq:fund_diag}, $\delta=\Dx$, and $T=5$ as the final time. We consider the network depicted in Fig.\ \ref{fig:network}, with $N_\mathcal R=8$, $N_\mathcal J=8$, and $N_\mathcal D=2$. 
We choose $\Dx=0.01$ and $\Dt=0.005$. The length of the roads are given in the figure as multiples of $\Dx$.
The destination node for the first group is $\dd_1=\jj_7$, for the second group is $\dd_2=\jj_8$. 
It is also convenient setting $\Delta\tau=\Dt$ so that the discretization in time of the functions $\Next_d$ requires no \emph{ad hoc} procedure.

At the initial time the network is empty. Boundary conditions are zero at every nodes but $\jj_1$ and $\jj_3$. At $\jj_1$ we impose $\rho_1=0.3$, so that only the group 1 is present. At $\jj_3$ instead we impose $\rho_2=0.4$, so that only the group 2 is present. Boundary conditions are active for all $t\in[0,T]$ in the case of basic and rational behavior, while in the case of the highly rational behavior the inflow stops at $t=1$. Observe that, due to the particular choice of the network and of the initial data, drivers of the group 2 have a unique path leading to their destination, that is $\textsc{r}_4 \cup \textsc{r}_6 \cup \textsc{r}_8$. On the other hand, drivers of group 1 can choose between $\textsc{r}_1 \cup \textsc{r}_3 \cup \textsc{r}_6 \cup \textsc{r}_7$ and $\textsc{r}_1 \cup \textsc{r}_2 \cup \textsc{r}_5$. This also means that $\Next_1$ can change in time at $\jj_2$, while $\Next_2$ is always constant at every junction of the network.

The road weights $w_\rr$ are all chosen to coincide with the travel time, i.e. $w_\rr=\Upsilon_\rr$.

In the next figures we show the two densities $\rho_1$ and $\rho_2$ separately on roads $\rr_1$, $\rr_2$, $\rr_3$, $\rr_4$, and $\rr_6$.

\subsection*{Basic behavior} 
\begin{figure}[t!]
\begin{center}
\includegraphics[width=1\textwidth, trim = 35mm 25mm 20mm 20mm,clip]{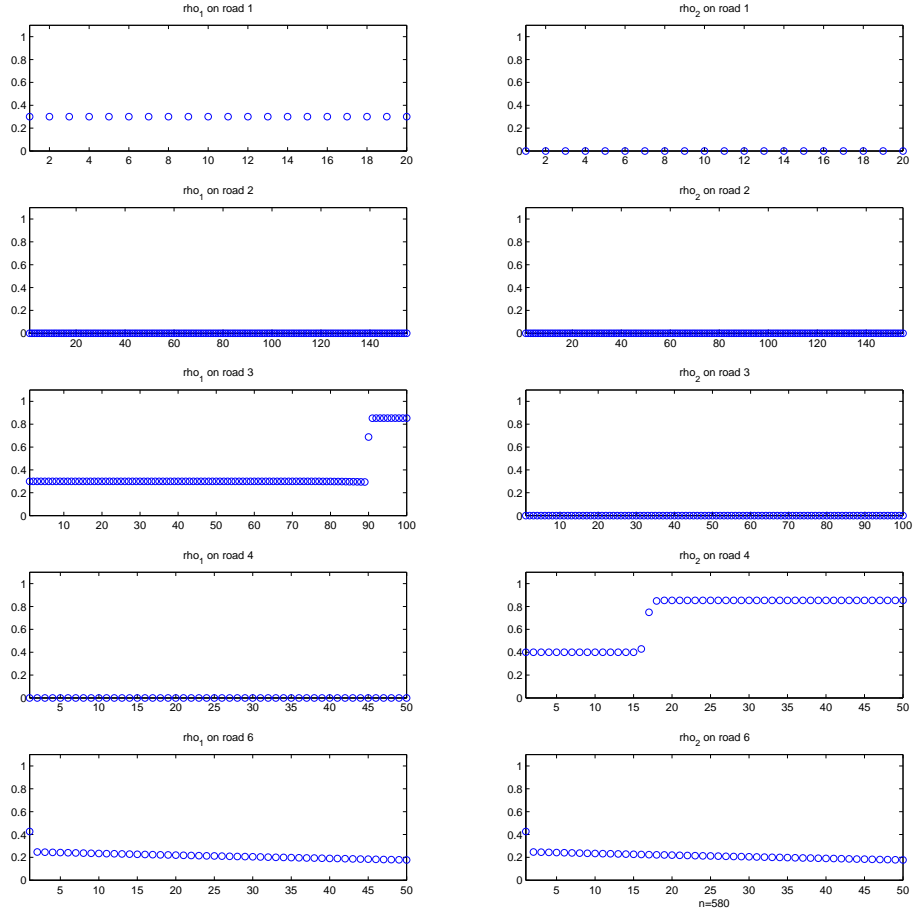}
\end{center}
\caption{Basic behavior.}
\label{fig:basic}
\end{figure}

In Fig.\ \ref{fig:basic} we show the result for the basic behavior at time $t=2.9$. As expected, this behavior prescribes drivers to follow the shortest path, which is, for group 1 entering at junction $\textsc{j}_1$, $\textsc{r}_1 \cup \textsc{r}_3 \cup \textsc{r}_6 \cup \textsc{r}_7$ and for group 2 entering at junction $\textsc{j}_3$, $\textsc{r}_4 \cup \textsc{r}_6 \cup \textsc{r}_8$. Then, roads $\textsc{r}_2$ and $\textsc{r}_5$ are unused. Note that the road $\textsc{r}_6$ is not able to gather the flows coming from roads $\textsc{r}_3$ and $\textsc{r}_4$, when drivers of both groups arrive. When this happens, two queues are formed along these incoming roads (you can see them in the figure). Note that the queues have the same level of density, as for the properties of the multi-path scheme \eqref{schema}, but travels backward at different velocities, due to the different boundary condition. Correctly, at the beginning of road $\rr_6$ we have the maximal flux, corresponding to the density $\rho_1+\rho_2=0.25+0.25=0.5$.

\subsection*{Rational behavior} 
\begin{figure}[t!]
\begin{center}
\includegraphics[width=1\textwidth, trim = 35mm 25mm 20mm 20mm,clip]{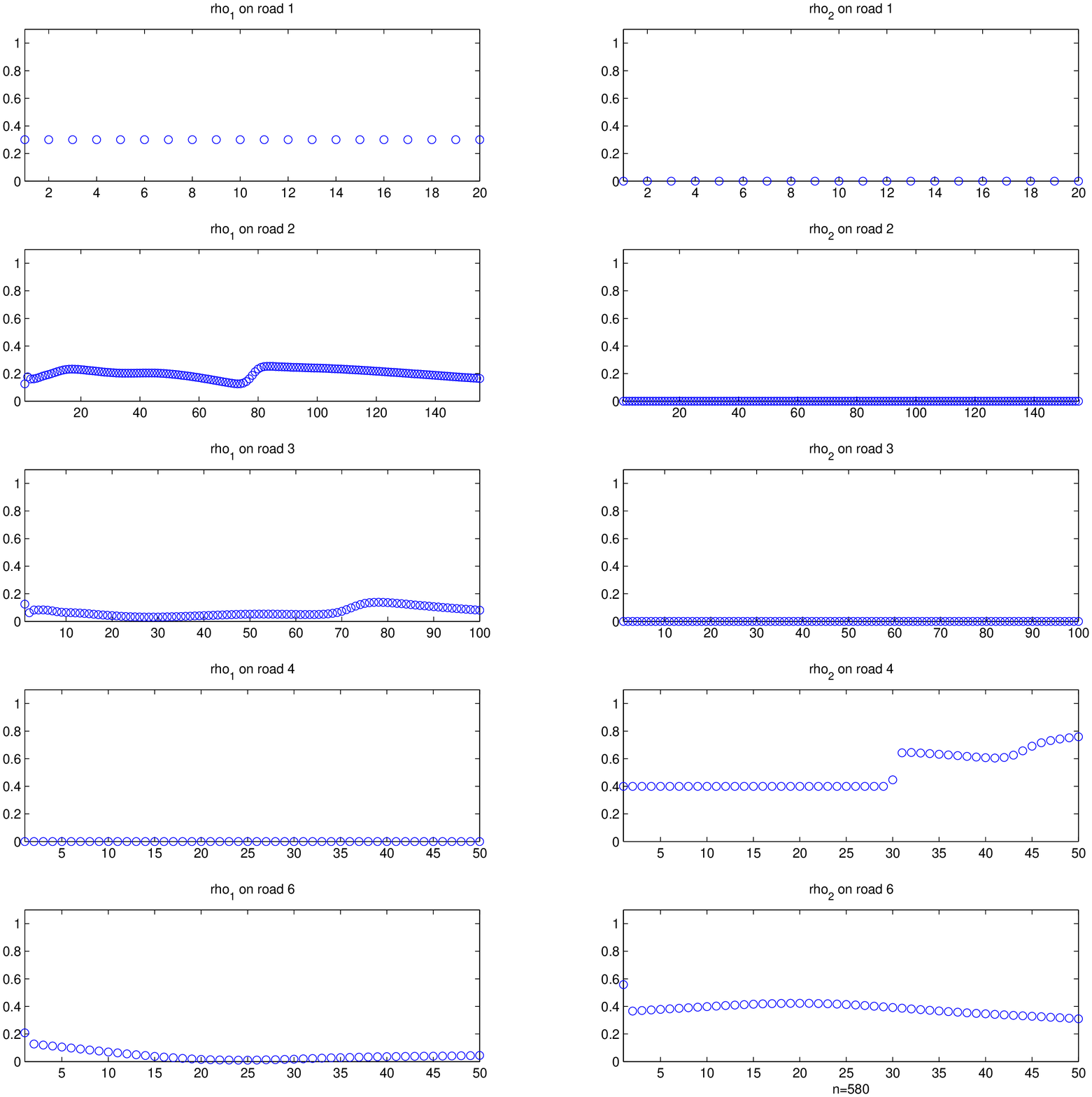}
\end{center}
\caption{Rational behavior.}
\label{fig:rational}
\end{figure}

In Fig.\ \ref{fig:rational} we show the result for the rational behavior at time $t=2.9$. The dynamics is more interesting than in the previous case. At first, drivers of  group 1 entering at junction $\textsc{j}_1$ will choose road $\textsc{r}_3$ because the path is shorter, but after some time the cars present in roads $\textsc{r}_3$ and $\textsc{r}_6$ (also including those belonging to group 2 which arrive from road $\textsc{r}_4$) make this path as ``expensive'' as the one passing in roads $\textsc{r}_2$ and $\textsc{r}_5$. Therefore, the drivers' preference at $\jj_2$ starts to oscillate between the two possibilities, generating in roads $\textsc{r}_2$ and $\textsc{r}_3$ the waves with positive speed that can be seen in the figure. The dynamics of group 2 is the same as before, but we observe that, thanks to the reduced number of cars of group 1 coming from $\textsc{r}_3$, the queue forming at junction $\textsc{j}_5$ and propagating in $\textsc{r}_4$ has a smaller level of density than in the previous case. Finally, no queue is formed in road $\textsc{r}_3$.

\subsection*{Highly rational behavior} 
\begin{figure}[ht!]
\begin{center}
\includegraphics[width=1\textwidth, trim = 35mm 25mm 20mm 20mm,clip]{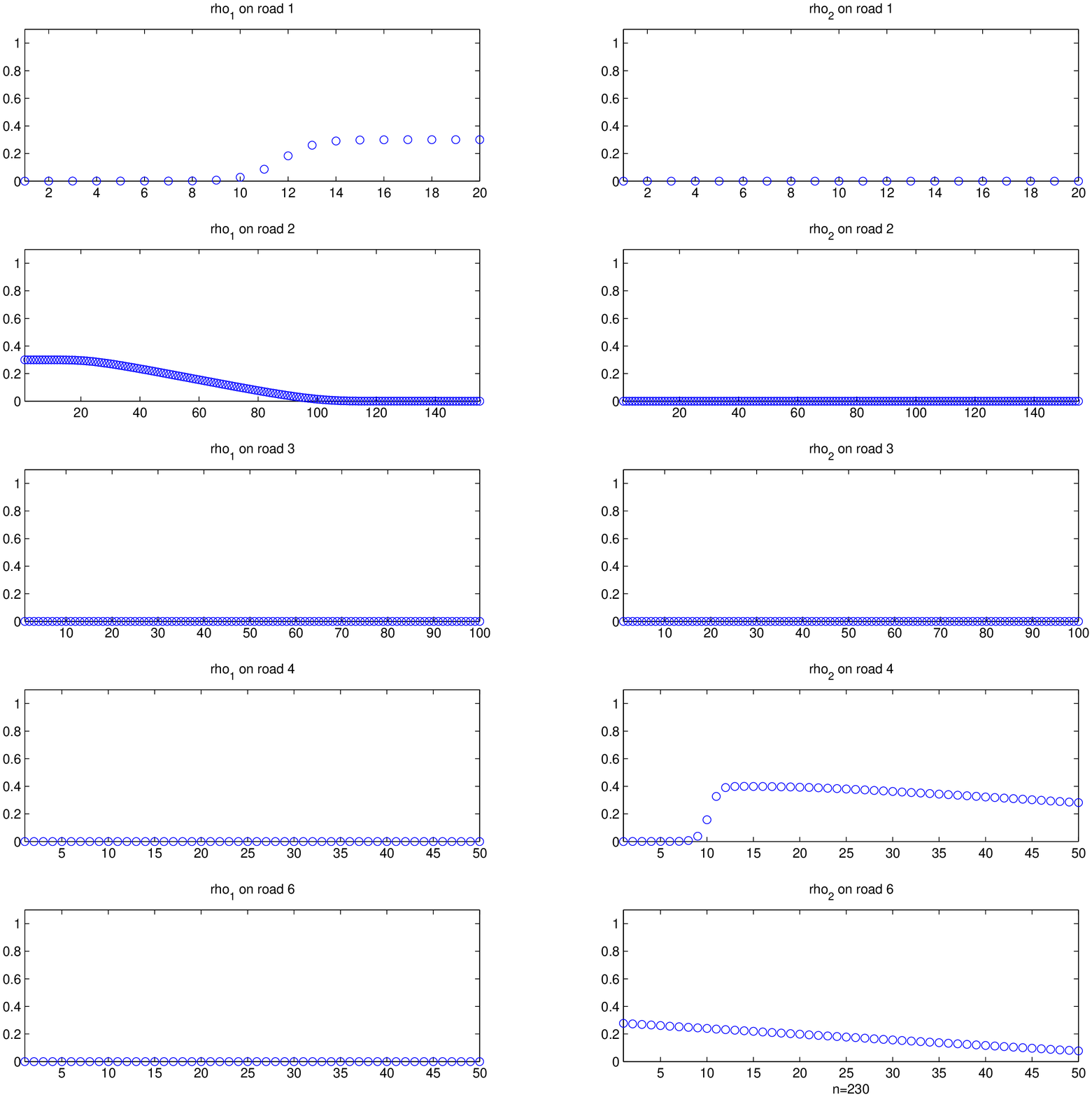}
\end{center}
\caption{Highly rational behavior, case A.}
\label{fig:hrational-A}
\end{figure}

In this case the algorithm oscillates between two solutions. The first one is shown in Fig.\ \ref{fig:hrational-A}, where we depicted the situation corresponding to the time $t=1.15$. Here, the drivers of group 1 entering at junction $\textsc{j}_1$ all choose to use road $\textsc{r}_2$, because they can forecast that a group of drivers will occupy $\textsc{r}_6$ at a later time (those of group 2 coming from $\textsc{r}_4$) and that such appearance will make the shorter path less convenient. Drivers of group 2 will follow their usual path $\textsc{r}_4 \cup \textsc{r}_6 \cup \textsc{r}_8$, but they will be undisturbed and no queue will appear at any junction. The jump shown in the figure for $\rho_2$ along $\textsc{r}_4$ is due to the boundary condition, because at $t=1$ the inflow for $\rho_2$ has passed from 0.4 to 0.

\begin{figure}[ht!]
\begin{center}
\includegraphics[width=1\textwidth, trim = 35mm 25mm 20mm 20mm,clip]{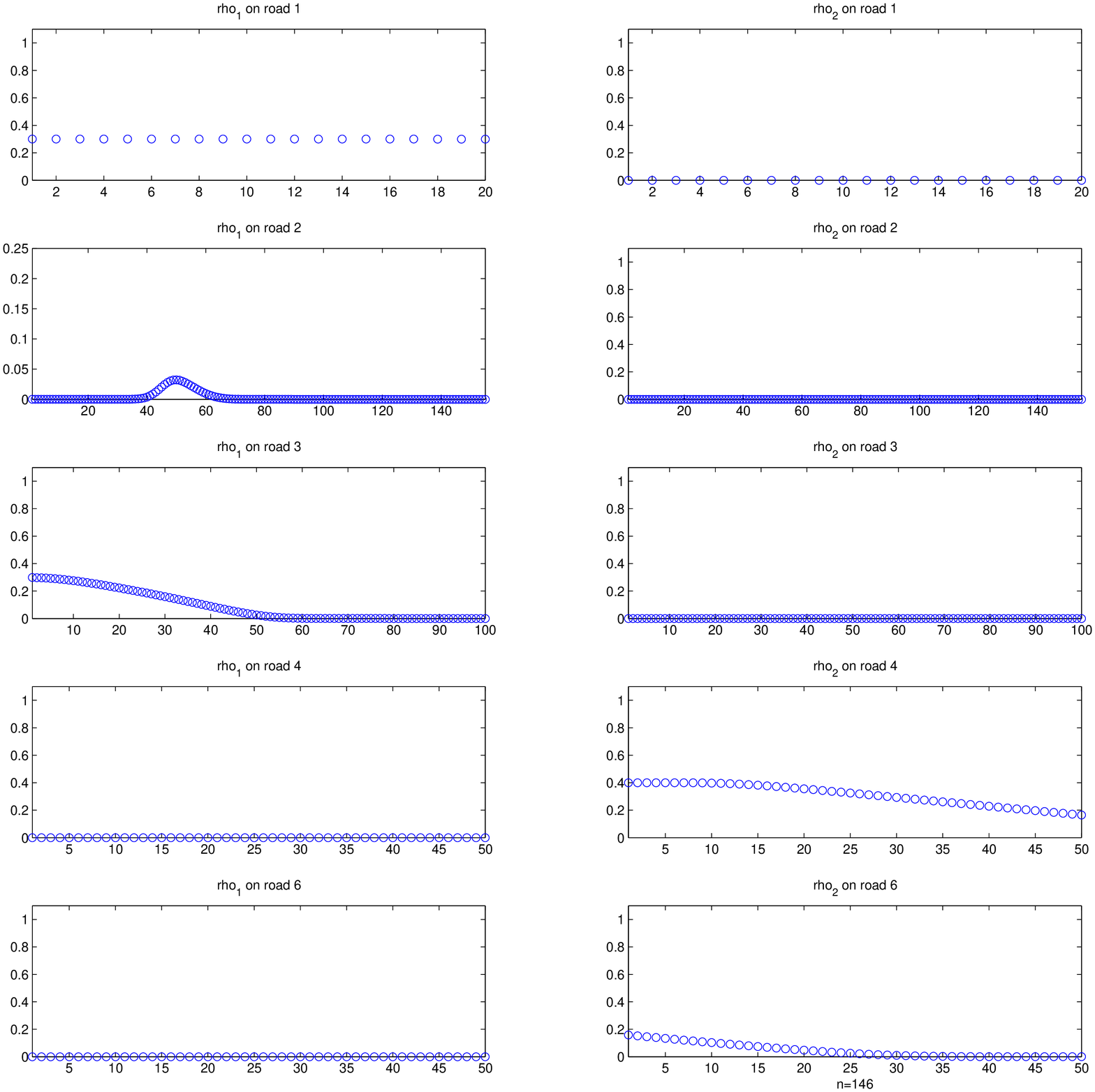}
\end{center}
\caption{Highly rational behavior, case B.}
\label{fig:hrational-B}
\end{figure}

The second solution is shown in Fig.\ \ref{fig:hrational-B}, where we depicted the situation corresponding to the time $t=0.73$. In this case, drivers of group 1 entering at junction $\textsc{j}_1$ will first choose $\textsc{r}_2$ (generating the wave that can be seen in the corresponding figure) and then they switch to $\textsc{r}_3$. The choice of switching time is performed so that drivers arrive at $\textsc{j}_5$ when the majority of drivers of group 2 has already passed the junction and is traveling along road $\textsc{r}_6$, reducing their effect on the dynamics for $\rho_1$. At a later stage a small queue will form in road $\textsc{r}_4$ when drivers of group 1 arrive to junction $\textsc{j}_5$, because the merging of the two groups exceeds the flow capacity of $\textsc{r}_6$. However, such queue disappears almost immediately because there are only a few cars of the second group remaining in $\textsc{r}_4$ (recall that the inflow ceases at $t=1$) and they soon pass through the junction.
%
%
%
%
%
%
\section{Conclusions and open problems}\label{sec:conclusions}
In this paper we have presented a novel model for traffic flows on a network, where drivers can modify their path during the evolution while preserving their initial destination. The model consists in solving separately the problem close to junctions and the one away from them, and then in suitably matching the solutions at the interfaces. 
The theoretical properties of the resulting combined model will be studied more in depth in a forthcoming paper, but the numerical results seem to be encouraging. A point of particular interest is whether the instability phenomena pointed out in~\cite{BressanYu} can occur also in our model. This is unclear at the moment because, as described in detail in~\cite{briani2014sub}, the adjacent cell to each junction \emph{behaves like a discretization of a buffer} and therefore its presence might have the same effect as the \emph{single buffer junctions} described in~\cite{BressanKhai} to obtain a well-posed problem on the network $\mathcal{N}$. In our view, however, the model is interesting by itself, being to our knowledge the first model allowing for changes in the path of groups of drivers during the evolution. 

The second novelty of this work is the introduction of different procedures to describe the choices of drivers at junctions, so to describe different levels of rationality (basic, rational and highly rational) as it was done in~\cite{cristiani2014sub} for pedestrian flows.
In this context it is still unclear under which assumptions the highly rational model, which is of particular interest since it leads to a Wardrop equilibrium on the network, admits a solution.
Another interesting open problem regards the stability of Wardrop equilibria. Numerical tests have shown that in some networks the algorithm oscillates between two non-equilibrium solutions, even if a small perturbation of an equilibrium solution is used as initial guess. 
On the other hand, we never observed oscillations among three or more solutions, even considering more complex networks with multiple complex junctions (not limited to the case $1\times 2$ and $2\times 1$ depicted in Figure \ref{fig:network}). 
A broader theoretical understanding about the relationship between the network and the equilibria seems to be necessary to proceed in this research field.  
Even without this knowledge, however, we believe that the introduction of rationality performed here can be a first step in a more ambitious program: to improve network design and network control so as to steer drivers' choices towards higher rationality than the amount of information they have available could allow. 
Following~\cite{cristiani2014sub}, we would like to adopt a network-scale control procedure to reproduce a rational evolution, even in regimes of reduced rationality. This would allow, for instance, to avoid undesired phenomena like Braess paradox (cf.~\cite{colombo2014}).

Similar considerations could be adapted also to models describing supply chains or multi-commodity flow problems. In these situations the control problems which give different degrees of rationality can describe some control procedure operated at junctions of the network by human beings who e.g. have to choose the destination of the various goods, based on different amounts of information about the structure and the current usage of each arc. Network optimization instead can be used by the network manager so as to pilot the choices of the operators at junctions towards a desirable global distribution, by acting on the capacity of the various arcs.


\medskip
Received xxxx 20xx; revised xxxx 20xx.
\medskip


\begin{thebibliography}{99}

\bibitem{colombosbg} 
S. Benzoni-Gavage, R. M. Colombo, 
\emph{An $n$-populations model for traffic flow},
Euro. J. Appl. Math., {\bf 14} (2003), 587--612.

\bibitem{BressanHan1} 
A. Bressan, K. Han, 
\emph{Optima and equilibria for a model of traffic flow},
SIAM J. Math. Anal., {\bf 43} (2011), 2384--2417.

\bibitem{BressanHan2} 
A. Bressan, K. Han, 
\emph{Nash equilibria for a model of traffic flow with several groups of drivers},
ESAIM Control. Optim. Calc. Var., {\bf 18} (2012), 969--986.

\bibitem{BressanHan3} 
A. Bressan, K. Han, 
\emph{Existence of optima and equilibria for traffic flow on networks},
Netw. Heterog. Media, {\bf 8} (2013), 627--648.

\bibitem{BressanKhai} 
A. Bressan, K. T. Nguyen, 
\emph{Conservation law models for traffic flow on a network of roads}, 
preprint (online version: {\tt http://www.math.ntnu.no/conservation/2014/009.pdf}).

\bibitem{BrPr} 
A. Bressan, F. S. Priuli, 
\emph{Infinite horizon noncooperative differential games}, 
J. Differential Equations, {\bf 227} (2006), 230--257.

\bibitem{BressanYu} 
A. Bressan, F. Yu, 
\emph{Continuous Riemann solvers for traffic flow at a junction}, 
preprint (online version: {\tt http://www.math.ntnu.no/conservation/2014/008.pdf}).

\bibitem{bretti2014DCDS-S}
G. Bretti, M. Briani, E. Cristiani, 
\emph{An easy-to-use algorithm for simulating traffic flow on networks: Numerical experiments}, 
Discrete Contin. Dyn. Syst. Ser. S, \textbf{7} (2014), 379--394. 

\bibitem{briani2014sub}
M. Briani, E. Cristiani, 
\emph{An easy-to-use numerical algorithm for simulating traffic flow on networks: Theoretical study},
Netw. Heterog. Media, \textbf{9} (2014), 519--552.

\bibitem{CaCr} 
S. Cacace, E. Cristiani, M. Falcone, 
\emph{Numerical approximation of Nash equilibria for a class of non-cooperative differential games}, 
In: L. Petrosjan, V. Mazalov (eds.), Game Theory and Applications, Vol. 16, Chap. 4, Nova Publishers, New York, 2013.

\bibitem{carlier2008SIAM}
G. Carlier, C. Jimenez, F. Santambrogio, 
\emph{Optimal transportation with traffic congestion and Wardrop equilibria},
SIAM J. Cont. Opt., {\bf 47} (2008), 1330--1350.

\bibitem{carlier2012JMS}
G. Carlier, F. Santambrogio, 
\emph{A continuous theory of traffic congestion and Wardrop equilibria},
J. Math. Sci., {\bf 181} (2012), 792--804.

\bibitem{cascone2007M3AS}
A. Cascone, C. D'Apice, B. Piccoli, L. Rarit\`a, 
\emph{Optimization of traffic on road networks}, 
Math. Models Methods Appl. Sci., \textbf{17} (2007), 1587--1617.

\bibitem{colombo2014}
R. M. Colombo, H. Holden,
\emph{On the Braess paradox with nonlinear dynamics and control theory}, 
preprint (online version: {\tt http://www.math.ntnu.no/conservation/2014/012.pdf}).

\bibitem{Cong}
Z. Cong, B. De Schutter, R. Babu\v{s}ka, 
\emph{Ant colony routing algorithm for freeway networks}, 
Transport. Res. C, {\bf 37C} (2013), 1--19.

\bibitem{cristiani2014sub}
E. Cristiani, F. S. Priuli, A. Tosin, 
\emph{Modeling rationality to control self-organization of crowds: An environmental approach}, 
SIAM J. Appl. Math., to appear.

\bibitem{cutolo2011}
A. Cutolo, C. D'Apice, R. Manzo, 
\emph{Traffic optimization at junctions to improve vehicular flows}, 
International Scholarly Research Network ISRN Applied Mathematics 01/2011, 2011.

\bibitem{dogbe2010MCM}
C. Dogb\'{e}, 
\emph{Modeling crowd dynamics by the mean-field limit approach}, 
Math. Comput. Modelling, {\bf 52} (2010), 1506--1520.

\bibitem{Fisk}
C. S. Fisk, 
\emph{Game theory and transportation systems modelling}, 
Transport. Res. B, {\bf 18B} (1984), 301--313.

\bibitem{herty2006SIOPT}
A. F\"ugenschuh, M. Herty, A. Klar, A. Martin, 
\emph{Combinatorial and continuous models for the optimization of traffic flows on networks}, 
SIAM J. Optim., {\bf 16} (2006), 1155--1176.

\bibitem{GaravelloDCDS}
M. Garavello, 
\emph{The LWR traffic model at a junction with multibuffers}, 
Discrete Contin. Dyn. Syst. S, {\bf 7} (2014), 463--482.

\bibitem{GaravelloGoatin}
M. Garavello, P. Goatin, 
\emph{The Cauchy problem at a node with buffer},
Discrete Contin. Dyn. Syst., {\bf 32} (2012), 1915--1938.

\bibitem{garavello2005CMS}
M. Garavello, B. Piccoli,
\emph{Source-destination flow on a road network}, 
Comm. Math. Sci., {\bf 3} (2005), 261--283.

\bibitem{piccolibook}
M. Garavello, B. Piccoli,
Traffic flow on networks, 
{\it AIMS Series on Applied Mathematics}, 
Springfield, MO, 2006.

\bibitem{gugat2005JOTA}
M. Gugat, M. Herty, A. Klar, G. Leugering,
\emph{Optimal control for traffic flow networks},
J. Optim. Theory Appl., {\bf 126} (2005), 589--616.

\bibitem{herty2003SISC}
M. Herty, A. Klar, 
\emph{Modeling, simulation, and optimization of traffic flow networks}, 
SIAM J. Sci. Comput., {\bf 25} (2003), 1066--1087.

\bibitem{herty2009NHM}
M. Herty, J.-P. Lebacque, S. Moutari,
\emph{A novel model for intersections of vehicular traffic flow},
Netw. Heterog. Media, {\bf 4} (2009), 813--826.

\bibitem{Hollander} 
Y. Hollander, J. N. Prashker, 
\emph{The applicability of non-cooperative game theory in transport analysis}, 
Transportation, {\bf 33} (2006), 481--496.

\bibitem{LachaWolf} 
A. Lachapelle, M.-T. Wolfram, 
\emph{On a mean field game approach modeling congestion and aversion in pedestrian crowds}, 
Transportation Res. B, {\bf 45} (2011), 1572--1589.

\bibitem{LW} 
M. J. Lighthill, G. B. Whitham, 
\emph{On kinematic waves II. A theory of traffic flow on long crowded roads}, 
Proc. Roy. Soc. London Ser. A, {\bf 229} (1955), 317--345.

\bibitem{Nachtigall} 
K. Nachtigall, 
\emph{Time depending shortest-path problems with applications  to railway networks},
Euro. J. Oper. Res., {\bf 83} (1995), 154--166.

\bibitem{OR1990} 
A. Orda, R. Rom, 
\emph{Shortest-path and minimum-delay algorithms in networks with time-dependent edge-length}, Journal of ACM, {\bf 37} (1990), 607--625.

\bibitem{Pr} 
F. S. Priuli, 
\emph{Infinite horizon noncooperative differential games with non-smooth costs}, 
J. Math. Anal. Appl., {\bf 336} (2007), 156--170.

\bibitem{priuli2014sub} 
F. S. Priuli, 
\emph{First order mean field games in pedestrian dynamics}, 
submitted. arXiv:1402.7296.

\bibitem{Richards} 
P. I. Richards, 
\emph{Shock waves on the highway}, 
Operations Res., {\bf 4} (1956), 42--51.

\bibitem{Wardrop} 
J. G. Wardrop, 
\emph{Some theoretical aspects of road traffic research}, 
Proc. Inst. Civ. Eng. Part II, {\bf 1} (1952), 325--378.


\end{thebibliography}
\end{document}